\documentclass[preprint,review,3p,times]{elsarticle}

\usepackage[mode=multiuser,draft]{fixme} 
\FXRegisterAuthor{cc}{envcc}{CC}
\FXRegisterAuthor{gk}{envgk}{GK}
\fxusetheme{color}

\usepackage{geometry}

\usepackage{tikz,tikzscale}
\usepackage{graphicx}

\usetikzlibrary{calc}
\usetikzlibrary{shadings}
\usetikzlibrary{shapes.geometric}
\usetikzlibrary{shadows.blur}
\usetikzlibrary{positioning}
\usetikzlibrary{backgrounds}
\usetikzlibrary{matrix}
\usetikzlibrary{intersections}

\usepackage{pgfplots}
\usepgfplotslibrary{groupplots}
\usepgfplotslibrary{fillbetween}
\usepackage{pgfplotstable}
\pgfplotsset{compat=1.9}
\pgfplotsset{every tick label/.append style={font=\footnotesize}}

\pgfplotsset{p1/.style={orange,mark=o, semithick}}
\pgfplotsset{p2/.style={green,mark=otimes, semithick}}
\pgfplotsset{p3/.style={red,mark=*, semithick}}
\pgfplotsset{p4/.style={blue,mark=oplus, semithick}}

\pgfplotsset{p11/.style={orange,mark=o, semithick, dashed}}
\pgfplotsset{p22/.style={green,mark=otimes, semithick, dashed}}
\pgfplotsset{p33/.style={red,mark=*, semithick, dashed}}
\pgfplotsset{p44/.style={blue,mark=oplus, semithick, dashed}}

\pgfplotsset{padding/.style={orange,mark=otimes, mark size=2pt}}
\pgfplotsset{basic/.style={black,mark=o, mark size=2pt}}
\pgfplotsset{cache/.style={blue!60!,mark=*, mark size=2pt}}
\pgfplotsset{cf/.style={red!60!,mark=square*, mark size=2pt}}

\usepackage{algorithm}
\usepackage{algpseudocodex}

\usepackage{xspace}
\usepackage{bold-extra}
\usepackage[most]{tcolorbox}

\usepackage{amsmath}
\usepackage{amsfonts}

\usepackage{booktabs}

\usepackage{caption}
\usepackage[labelformat=simple]{subcaption}

\usepackage{stmaryrd}
\SetSymbolFont{stmry}{bold}{U}{stmry}{m}{n}

\usepackage{braket,amsfonts}

\usepackage{array}

\usepackage{mathtools}

\usepackage{arydshln}

\usepackage{multirow}

\usepackage{amssymb}

\journal{Journal of Computational Physics}

\begin{document}

\begin{frontmatter}



\title{Multigrid methods for the Stokes problem on GPU systems} 


\author[a]{Cu Cui\corref{cor1}} 
\ead{cu.cui@iwr.uni-heidelberg.de}
\author[a]{Guido Kanschat} 

\cortext[cor1]{Corresponding author}

\affiliation[a]{
organization={Interdisciplinary Center for Scientific Computing (IWR), Heidelberg University},
            addressline={Im Neuenheimer Feld 205}, 
            city={Heidelberg},
            postcode={69120}, 
            country={Germany}}

\begin{abstract}
This paper presents a matrix-free multigrid method for solving the Stokes problem, discretized using $H^{\text{div}}$-conforming discontinuous Galerkin methods. We employ a Schur complement method combined with the fast diagonalization method for the efficient evaluation of the local solver within the multiplicative Schwarz smoother. This approach operates directly on both the velocity and pressure spaces, eliminating the need for a global Schur complement approximation. By leveraging the tensor product structure of Raviart-Thomas elements and an optimized, conflict-free shared memory access pattern, the matrix-free operator evaluation demonstrates excellent performance numbers, reaching over one billion degrees of freedom per second on a single NVIDIA A100 GPU. Numerical results indicate efficiency comparable to that of the three-dimensional Poisson problem.
\end{abstract}




\end{frontmatter}

\section{Introduction}\label{sec:intro}

The efficient solution of the Stokes equations is a crucial task in computational fluid dynamics (CFD), especially in stationary computations or applying implicit-explicit timestepping schemes where the inertia term of the Navier-Stokes equations is treated explicitly.
Multigrid methods are known to be among the most efficient preconditioners and solvers. Combining divergence-conforming discontinuous Galerkin (DG) methods with overlapping vertex-patch smoothers yields effective multigrid methods and achieves convergence rates for the Stokes problem which are comparable to those for the Laplacian~\cite{kanschat2015multigrid}. 
In recent years, the advent of Graphics Processing Units (GPUs) has provided significant computational power, but algorithms have to be redesigned to make use of these devices. Hence, in this work, we focus on the GPU acceleration such a multigrid method. This work can be seen as part of a global effort to provide efficient GPU implementations entering popular computational libraries such as libCEED~\cite{Abdelfattah2021}, PETSc~\cite{mills2021toward}, MFEM~\cite{anderson2021mfem} and deal.II~\cite{dealII95}. 

The robustness and efficiency of DG methods for incompressible flows and their implementation on CPU is well studied~\cite{krank2017high,fehn2017stability,fehn2018efficiency}.
The GPU-based incompressible fluid flow simulator NekRS~\cite{Kolev2021,Merzari2023} has achieve unprecedented simulation scales and accuracy on high-performance computing platforms with backward difference temporal discretization and high-order spectral elements for spatial discretization, a $p$-multigrid method with Chebyshev-accelerated additive Schwarz smoother is employed for solving the pressure-Poisson problem.

 For GPUs, the optimized implementation of first-order hyperbolic systems is well-documented. 
The works~\cite{franco2020high,karakus2019gpu} considered multigrid methods with Chebyshev or ILU smoothers for the Laplacian and demonstrated close to empirically predicted roofline performance, where the matrix-free~\cite{KronbichlerKormann19} implementation of these operators is often crucial for achieving high performance. We also refer to our work on the GPU implementation of vertex-patch smoothers for the Laplacian~\cite{Cui2023,Cui2024,cui2024acceleration}.
However, an adaptation of the Stokes operator remains relatively unexplored.

In~\cite{Will2023}, a refined low-order method was introduced for preconditioning grad-div problems within the $H^{\text{div}}$ space. The study demonstrated that the high-order operator, discretized using a so-called interpolation–histopolation basis, is spectrally equivalent to the lowest-order discretization applied on a refined mesh. Subsequent work~\cite{Will2024} expanded on this by developing efficient, matrix-free preconditioners for high-order $H^{\text{div}}$ problems, focusing on the associated $2 \times 2$ saddle-point system. These saddle-point systems were efficiently preconditioned using block preconditioners with GPU acceleration.

For the rather complex and expensive patch smoother, we demonstrated its efficient implementation on GPUs in our previous work on the Poisson equation~\cite{Cui2023,Cui2024}. In this study, we employ the Schur complement method combined with fast diagonalization method~\cite{lynch1964direct} for the efficient evaluation of the local solver. Bank conflicts are one of the most significant issues in GPU computing, particularly when using shared memory. We have adapted methods from our previous work to design a conflict-free memory access pattern, which significantly enhances performance, achieved over one billion degrees of freedom.

 In the present work, we propose a matrix-free method for the Stokes operator, leveraging Raviart–Thomas elements and employing detailed GPU optimization strategies. This approach operates on both velocity and pressure spaces, eliminating the need for a global Schur complement approximation.
 
The remainder of this paper is organized as follows. In Section~\ref{sec:stokes}, we present the mathematical formulation for the DG discretization of the approximate Stokes equations. Section~\ref{sec:MG} provides a detailed explanation of the multigrid method and the vertex patch smoother. This is followed by the numerical validation test cases discussed in Section~\ref{sec:tests}. In Section~\ref{sec:gpu}, we delve into the key aspects of the GPU implementation, performance analysis, and core kernel optimizations. Finally, Section~\ref{sec:con} concludes the paper.

\section{The Stokes problem}\label{sec:stokes}

We examine the numerical solutions for the Stokes equations under the following formulation:
\begin{gather}\label{eq:model}
\begin{aligned}
    -\Delta \boldsymbol{u} + \nabla p = \boldsymbol{f} \qquad & \text{ in } \Omega \\
    \nabla \cdot \boldsymbol{u} = 0 \qquad & \text{ in } \Omega \\
                \boldsymbol{u} = 0 \qquad & \text{ on } \partial\Omega
\end{aligned}
\end{gather}
where the domain $\Omega \subset \mathbb{R}^d$ is bounded and convex, with dimensions $d = 2$ or $3$. We enforce no-slip boundary conditions on $\partial\Omega$. Suitable function spaces for this problem are $V = H^1_0(\Omega; \mathbb{R}^d)$ for the velocity vector $\boldsymbol{u}$, and $Q = L^2_0(\Omega)$, which consists of square integrable functions with zero mean, for the pressure $p$. It is noteworthy that this setting is flexible and can accommodate other types of well-defined boundary conditions without compromising the well-posedness of the problem.

\subsection{Discontinuous Galerkin discretization}\label{sec:DG}
To implement a finite element discretization of the Stokes equations, we first subdivide the domain $\Omega$ into a mesh $\mathcal{T}_h$, consisting of parallelogram cells in 2D and parallelepiped cells in 3D. We define the combined discrete spaces as $X = V \times Q$. In alignment with the methodology discussed in~\cite{cockburn2007note}, we select $V$ from the subspace $H^{\text{div}}_0(\Omega)$, characterized as follows:
\begin{gather*}
H^{\text{div}}(\Omega) = \{v\in L^2(\Omega;\mathbb{R}^d) | \nabla\cdot v \in L^2(\Omega)\}, \\
H^{\text{div}}_0(\Omega) = \{v\in H^{\text{div}}(\Omega) | v\cdot\boldsymbol{n}=0 \text{ on }\partial\Omega\}.
\end{gather*}
For this implementation, we utilize the Raviart–Thomas space~\cite{Raviart1997}, although other combinations of velocity and pressure spaces that satisfy the relation $\nabla \cdot V = Q$ are also valid. 

Within reference cell, we construct polynomial spaces: $\widehat{Q}_k$, consisting of polynomials of degree up to $k$ in $d$ dimensions, and $\widehat{V}_k$, the space of vector-valued Raviart–Thomas polynomials, defined as $\widehat{Q}^d_k + x\widehat{Q}_k$. Here, the polynomial degree $k$ is uniformly applied across the mesh and henceforth we will drop the subscript $k$. The resulting finite element spaces are defined as:
\begin{gather*}
V_h = \{v\in H^{\text{div}}_0(\Omega) | \forall K\in\mathcal{T}_h : v_{|K} \in V_K \}, \\
Q_h = \{v\in L^{2}_0(\Omega) | \forall K\in\mathcal{T}_h : q_{|K} \in Q_K\}.
\end{gather*}
We denote Raviart-Thomas elements of degree $k \geq 0$ as $\text{RT}_k$, and use the shorthand $\mathbb{Q}_{k+1,k} \times \mathbb{Q}_{k,k+1}$ in 2D and a similar notation in 3D for the tensor product spaces. For $k=0$, we obtain the \emph{marker and cell} scheme, a staggered finite difference method, according to~\cite{Kanschat08mac}.

Although the normal components of functions within $V_h$ are continuous across cell interfaces, their tangential components are not, leading to $V_h \not\subset H^1(\Omega; \mathbb{R}^d)$. We follow the example in for instance~\cite{cockburn2007note,Kanschat20081093,KANSCHAT20105933} and apply a DG formulation to the discretization of the elliptic operator, namely the interior penalty method~\cite{Arnold1982742}. Let $K_1$ and $K_2$ be two mesh cells with a joint face $F$, and let $u_1$ and $u_2$ be the traces of a function $u$ on $F$ from $K_1$ and $K_2$, respectively. On this face $F$, we introduce the averaging operator
\begin{equation}\label{eq:average}
    \{\!\!\{ \boldsymbol{u} \}\!\!\} = \frac{\boldsymbol{u}_1 + \boldsymbol{u}_2}{2},
\end{equation}
and $\{\!\!\{ \boldsymbol{u} \}\!\!\} = \boldsymbol{u}$ on the boundary. In this notation, the interior penalty bilinear form reads
\begin{gather}\label{eq:bilinear}
\begin{aligned}
    a_h(\boldsymbol{u},\boldsymbol{v}) = & (\nabla \boldsymbol{u}, \nabla \boldsymbol{v})_{\mathcal{T}_h} + 4\langle \gamma_e \{\!\!\{\boldsymbol{u}\otimes\boldsymbol{n} \}\!\!\}, \{\!\!\{ \boldsymbol{v}\otimes\boldsymbol{n} \}\!\!\} \rangle_{{\mathcal{F}^i_h}} \\ 
    & -2\langle\{\!\!\{\nabla\boldsymbol{u}\}\!\!\}, \{\!\!\{ \boldsymbol{v}\otimes\boldsymbol{n} \}\!\!\} \rangle_{{\mathcal{F}^i_h}}
    -2\langle\{\!\!\{\nabla\boldsymbol{v}\}\!\!\}, \{\!\!\{ \boldsymbol{u}\otimes\boldsymbol{n} \}\!\!\} \rangle_{{\mathcal{F}^i_h}} \\
    & + 2\langle \gamma_e \boldsymbol{u}, \boldsymbol{v}\rangle_{{\mathcal{F}^{\partial}_h}}
    -\langle \partial_n\boldsymbol{u}, \boldsymbol{v}\rangle_{{\mathcal{F}^{\partial}_h}}
    -\langle \boldsymbol{u}, \partial_n\boldsymbol{v}\rangle_{{\mathcal{F}^{\partial}_h}}.
\end{aligned}
\end{gather}
The bilinear form $(a, b)_{\mathcal{T}_h} = \sum_{K \in \mathcal{T}_h}\int_K ab \text{ dx}$ denotes volume integrals and $\langle a, b \rangle_{\mathcal{F}_h} = \sum_{F \in \mathcal{F}_h}\int_F ab \text{ ds}$ boundary integrals. Here, the operator `$\otimes$' denotes the dyadic product of two vectors.

The discrete weak formulation of~\eqref{eq:model} reads now: find $(\boldsymbol{u}_h, p_h) \in V_h \times Q_h$, such that for all test functions $\boldsymbol{v}_h \in V_h$ and $q_h \in Q_h$ there holds:
\begin{equation}\label{eq:weak_form}
    \mathcal{A}_h \equiv a_h(\boldsymbol{u}_h,\boldsymbol{v}_h) + (p_h, \nabla \cdot \boldsymbol{v}_h) - 
    (q_h, \nabla \cdot \boldsymbol{u}_h) = (f, \boldsymbol{v}_h).
\end{equation}

\subsection{Matrix-free evaluation for Stokes Operator}\label{sec:mf}
With in each mesh cell $K \in \mathcal{T}_h$, we have that the velocity $\boldsymbol{u}_h^{K}$ and pressure $p_h^{K}$ representations are expanded as follows: $\boldsymbol{u}_h^{K}=\sum_{i=1}^N \boldsymbol{\varphi}_i^{K}(\mathbf{x}) \boldsymbol{u}_i^{K}$ and $p_h^{K}=\sum_{i=1}^N {\psi}_i^{K}(\mathbf{x}) p_i^{K}$. Here $\boldsymbol{u}_i^{K}$ and $p_i^{K}$ are the coefficients derived from the variational principles and $\boldsymbol{\varphi}_i^{K}(\mathbf{x})$ and $ {\psi}_i^{K}(\mathbf{x})$ are the corresponding basis functions, defined as polynomials $\boldsymbol{\varphi}_i(\hat{\mathbf{x}})$ and $ {\psi}_i^{K}(\hat{\mathbf{x}})$ on the reference element $\hat K$ with coordinates $\hat{\mathbf{x}}$. These are then mapped to the real coordinates $\mathbf{x}$ through the transformation $\boldsymbol{F}_K(\hat{\mathbf{x}})$ as $\mathbf{x}=\boldsymbol{F}_K(\hat{\mathbf{x}})$, with the Jacobian $\mathcal{J}_{K}(\hat{\mathbf{x}})=$ $\frac{\mathrm{d} \boldsymbol{F}_K}{\mathrm{d} \mathbf{x}}$, which in the Cartesian case simplifies to diagonal and constant within each cell, namely $\mathcal{J} = \text{diag}(\frac{1}{h},\frac{1}{h})$.

\begin{algorithm}[tp]
\caption{Matrix-free evaluation of Stokes operator~\eqref{eq:bilinear}.}\label{alg:mf}
\begin{itemize}
    \item[i)] loop over cells $\mathcal{T}_h$ 
        \begin{itemize}
        \item[a)] gather local vector values $\boldsymbol{u}^K_i$ and $p^K_i$ on cell from global input vector $\boldsymbol{u}$ and $p$
        \item[b)] interpolate local vector gradients $\nabla \boldsymbol{u}^K$ and local vector values $p^K$ in quadrature points, $\nabla \boldsymbol{u}^K_h(\textbf{x}_q) = \sum_i \hat{\nabla}\hat{\boldsymbol{\varphi}}_i(\hat{\textbf{x}}_q)\text{det}(\mathcal{J}^{-T})\boldsymbol{u}^K_i$, $p^K_h(x_q) = \sum_i\hat{\psi}_i(\hat{x}_q)p^K_i$
        \item[c)] for each quadrature index $q$, prepare integrand at each quadrature point by computing
        $\boldsymbol{r}_q = \mathcal{J}^{-1}\nabla \boldsymbol{u}^K_h(\textbf{x}_q) w_q$,
         $\boldsymbol{s}_q = \nabla\cdot \boldsymbol{u}^K_h(\textbf{x}_q) w_q$, $t_q = p^K_h(x_q) * w_q$
        \item[d)] evaluate local integrals for all test functions $i$ by quadrature $\boldsymbol{y}^K_i = (\nabla \boldsymbol{u}^K_h, \nabla \boldsymbol{\varphi}_i)_K + (p^K_h, \nabla \cdot \boldsymbol{\varphi}_i)_K \approx \sum_q \mathcal{J}\hat{\nabla}\hat{\boldsymbol{\varphi}}_i(\hat{\textbf{x}}_q) \boldsymbol{r}_q + \sum_q \hat{\nabla}\cdot\hat{\boldsymbol{\varphi}}_i(\hat{\textbf{x}}_q) t_q$, 
        $z^K_i = (\psi_i, \nabla\cdot \boldsymbol{u}^K_h(\textbf{x}_q)) \approx \sum_q \hat{\psi}_i(\hat{x}_q) \boldsymbol{s}_q$
        \end{itemize}
    \item[ii)] loop over interior faces $\mathcal{F}^i_h$
        \begin{itemize}
            \item[a$^+$)] gather local vector values $\boldsymbol{u}^+_i$ from global input vector $\boldsymbol{u}$ associated with interior cell $K^+$
            \item[a$^-$)] gather local vector values $\boldsymbol{u}^-_i$ from global input vector $\boldsymbol{u}$ associated with interior cell $K^-$
            \item[b$^+$)] interpolate $\boldsymbol{u}^+$ in face quadrature points $\boldsymbol{u}^+_h(\textbf{x}_q)$, $\nabla \boldsymbol{u}^+_h(\textbf{x}_q)$
            \item[b$^-$)] interpolate $\boldsymbol{u}^-$ in face quadrature points $\boldsymbol{u}^-_h(\textbf{x}_q)$, $\nabla \boldsymbol{u}^-_h(\textbf{x}_q)$
            \item[c)] for each quadrature index $q$, compute the numerical flux contribution\\
            $\boldsymbol{s}_q=\gamma_e(\boldsymbol{u}^-_h(\textbf{x}_q) - \boldsymbol{u}^+_h(\textbf{x}_q)) - \frac{1}{2}(\nabla \boldsymbol{u}^-_h(\textbf{x}_q) + \nabla \boldsymbol{u}^+_h(\textbf{x}_q))\boldsymbol{n}^-$,
            $\boldsymbol{t}_q = -\frac{\gamma_e}{2}(\boldsymbol{u}^-_h(\textbf{x}_q) - \boldsymbol{u}^+_h(\textbf{x}_q))$
            \item[d$^+$)] evaluate local integrals by quadrature $\boldsymbol{y}^+_i \approx -\sum_q \mathcal{J}\hat{\boldsymbol{\varphi}}_i(\hat{\textbf{x}}_q)\boldsymbol{s}_q w_q - \sum_q \hat{\nabla}\hat{\boldsymbol{\varphi}}_i(\hat{\textbf{x}}_q)\boldsymbol{n}^-\boldsymbol{t}_q$
            \item[d$^-$)] evaluate local integrals by quadrature $\boldsymbol{y}^-_i \approx \sum_q \mathcal{J}\hat{\boldsymbol{\varphi}}_i(\hat{\textbf{x}}_q)\boldsymbol{s}_q w_q - \sum_q \hat{\nabla}\hat{\boldsymbol{\varphi}}_i(\hat{\textbf{x}}_q)\boldsymbol{n}^-\boldsymbol{t}_q$
            \item[e$^+$)] add local contribution $\boldsymbol{y}^+_i$ into the global result vector $\boldsymbol{y}$ associated with $K^+$
            \item[e$^-$)] add local contribution $\boldsymbol{y}^-_i$ into the global result vector $\boldsymbol{y}$ associated with $K^-$
        \end{itemize}
    \item[iii)] loop over boundary faces $\mathcal{F}^{\partial}_h$
        \begin{itemize}
            \item[a)] gather local vector values $\boldsymbol{u}^-_i$ from global input vector $\boldsymbol{u}$ associated with interior cell $K^-$
            \item[b)] interpolate $\boldsymbol{u}^-$ in face quadrature points $\boldsymbol{u}^-_h(\textbf{x}_q)$, $\nabla \boldsymbol{u}^-_h(\textbf{x}_q)$
            \item[c)] for each quadrature index $q$, compute the numerical flux contribution and multiply by integration factor
            \item[d)] evaluate local integrals $\boldsymbol{y}^-_i$ by quadrature similar to inner faces
            \item[e)] add local contribution $\boldsymbol{y}^-_i$ into the global result vector $\boldsymbol{y}$ associated with $K^-$
        \end{itemize}
\end{itemize}
\end{algorithm}

Similar to the scalar cases~\cite{KronbichlerKormann19,Cui2024}, the Stokes operator can also be evaluated in a matrix-free manner. The matrix-free computation of the integrals, which represent the product of the matrix $\mathcal{A}$ with a specified vector, is facilitated through iterative loops across all cells and faces that are part of the operator described in equation~\eqref{eq:bilinear}.
The process outlined in Algorithm~\ref{alg:mf} delineates the steps involved for the Stokes operator. This evaluation is methodically divided into three distinct phases: one for cell integrals, another for integrals over interior faces, and a third for integrals over boundary faces, corresponding to the respective sums in equation~\eqref{eq:bilinear}.

These iterations are conceptually divided into five key stages, namely: (a) gathering the solution values relevant to the current cell(s), (b) computing the values or gradients of the local solution at the designated quadrature points, (c) applying the necessary geometric transformations at these quadrature points, (d) multiplying these results by the test functions followed by summation via numerical quadrature, and (e) integrating these local intermediate results and distributing the local vector into the global vector.

\subsection{Sum Factorization}\label{sec:sf}

In the matrix-free approach outlined in Algorithm~\ref{alg:mf}, particularly steps (ii-iv)(b) and (ii-iv)(d), interpolating solution coefficients to quadrature points and performing summations at these points are pivotal. These steps highlight how every vector component within a cell contributes to the values at each quadrature point. By leveraging the tensor product structure of the basis functions and quadrature formulas, the computational complexity can be reduced to $O(d k^{d+1})$ using sum factorization. This technique was originally introduced in the spectral element community, as described in~\cite{orszag1979spectral}, and later has been extended to discontinuous Galerkin (DG) methods, as illustrated in~\cite{KronbichlerKormann19}.
\begin{figure}
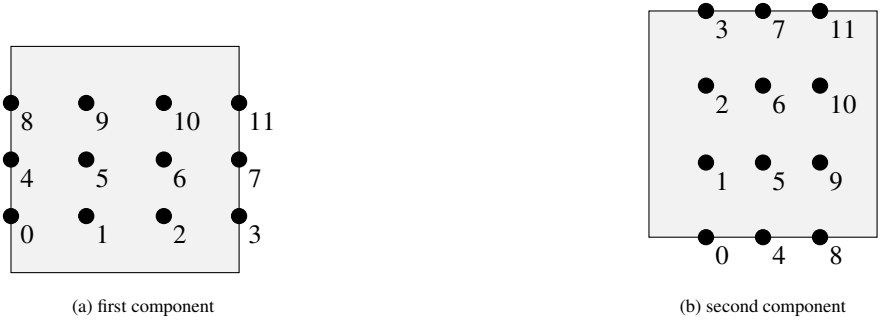

\begin{subfigure}[b]{0.49\textwidth}
\centering
    \includegraphics{Figures/rt_element.tikz}
    \caption{first component}
\end{subfigure}
\begin{subfigure}[b]{0.49\textwidth}
\centering
    \includegraphics{Figures/rt_element_y.tikz}
    \caption{second component}
\end{subfigure}
\caption{Visualization of degrees of freedom layout for RT$_2$ element in two dimensions. The direction in which the degrees of freedom are continuous is defined as the parallel direction, while the other direction (or directions in 3D) is referred to as the orthogonal direction.}
\label{fig:RT2}
\end{figure}
As depicted in Figure~\ref{fig:RT2}, Raviart-Thomas elements also exhibit a similar tensor product structure, allowing us to employ sum factorization to efficiently evaluate the Stokes operator~\cite{witte2022fast}.

It’s worth noting that optimizations are leveraged for cases where the local cell matrix—steps (b) through (d) in Algorithm~\ref{alg:mf}—can be expressed as a sum of Kronecker products, like the Laplacian in Cartesian geometries represented as $A = N_2 \otimes M_1 + M_2 \otimes N_1$ in 2D. This configuration simplifies sum-factorization operations, particularly under conditions of constant coefficients and axis-aligned meshes.

\section{Multigrid method}\label{sec:MG}

To formulate the finite element multilevel preconditioner as proposed by~\cite{braess1983,bramble1990parallel}, we establish a hierarchy of refined meshes:
\begin{gather*}
    \mathbb{T}_0 \sqsubset \mathbb{T}_1 \sqsubset\dots\sqsubset\mathbb{T}_L,
\end{gather*}
with each cell in $\mathbb{T}_\ell$ for $\ell \geq 1$ being a regular refinement of a cell from $\mathbb{T}_{\ell-1}$, resulting in the creation of $2^d$ child cells. In our experimental setup, $\mathbb{T}_0$ represents the refinement of a single cell into $2^d$ cells.
Associated with these subdivisions, we define a corresponding hierarchy of nested finite element spaces:
\begin{gather*}
    V_0 \subset V_1\subset \dots\subset V_L\\
    Q_0 \subset Q_1\subset \dots\subset Q_L\\
   X_0 = V_0 \times Q_0 \subset X_1\subset \dots\subset V_L \times Q_L = X_L,
\end{gather*}
Between the spaces $X_\ell$, we introduce the prolongation operator $I_{\ell}^{\uparrow}$ for transferring from space $X_{\ell}$ to $X_{\ell+1}$ via canonical embedding and the restriction operator $I_{\ell}^{\downarrow}$ for performing the $l_2$--projection from $X_{\ell+1}$ to $X_{\ell}$.
 
We consider the V-cycle as the most basic multigrid cycle and apply it to the Stokes system. The multigrid preconditioner $\mathcal{P}_\ell^{-1}$ is recursively defined starting from the coarsest level where $\mathcal{P}_0^{-1} = \mathcal{A}_0^{-1}$. Letting $x:=(\boldsymbol{u},p)^T$, the action of $\mathcal{P}^{-1}_\ell$ on a vector $b_\ell$ for level $\ell \geq 1$ proceeds as follows:
\begin{itemize}
    \item[(1)] Pre-smoothing: 
    \begin{equation*}
        x_\ell \gets S_\ell(x_\ell,b_\ell)
    \end{equation*} 
    \item[(2)] Coarse grid correction: 
    \begin{equation*}
        x_\ell \gets x_\ell + I_{\ell-1}^\uparrow \mathcal{P}^{-1}_{\ell-1} I_{\ell-1}^\downarrow (b_\ell - \mathcal{A}_\ell x_\ell)
    \end{equation*}
    \item[(3)] Post-smoothing:
    \begin{equation*}
        x_\ell \gets S_\ell(x_\ell,b_\ell)
    \end{equation*}
    \item[] and let $P^{-1}_{\ell-1}b_\ell = x_\ell$.
\end{itemize}
The smoothing operations $S_\ell$ on each level $\ell$ are essential for the algorithm’s efficiency, and specifics of these operations will be detailed later in the discussion. We consistently adopt a single pre-smoothing and post-smoothing step throughout this article. We note that certain theoretical frameworks require the variable V-cycle, but experimentally, our choice performs best.

\subsection{Vertex-patch smoothers}\label{sec:smoother}

\begin{figure}[tp]
\centering
\includegraphics[width=.4\textwidth]{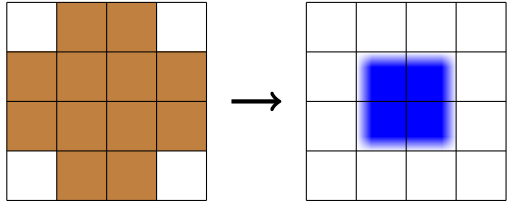}
\caption{Left: domain of dependence $\overline X_j$ (brown). Right: the support of $X_j$ (blue) of the vertex patch solver where fading indicates that degrees on the boundary of the patch are not part of $X_j$.}
\label{fig:smoother}
\end{figure}

In this section, we elaborate on a class of smoothing operators, denoted as $S_\ell$ which leverage a subspace decomposition of the space $X_\ell$. Let $\mathbb V_\ell$ represent the set of all interior vertices within the mesh $\mathbb{T}_\ell$ at level $\ell$ of the multilevel discretization. Each vertex $p_j \in \mathbb V_\ell$ corresponds to a vertex patch comprising all cells that share $p_j$ as a vertex. For structured quadrilateral and hexahedral meshes, these vertex patches typically encompass 4 and 8 cells, respectively. Collectively, they provide an overlapping coverage of the entire mesh $\mathbb{T}_\ell$. The total number of these patches is denoted by $J$, which is the cardinality of $\mathbb V_\ell$.

Each vertex patch $j$ is associated with a subspace $X_{\ell,j}$ of the finite element solution space, consisting of functions supported exclusively on that patch and zero elsewhere. Note that this configuration naturally enforces no-slip boundary conditions of velocity subspace $V_{\ell,j}$. The local solvers $P_{\ell,j} : X_\ell \rightarrow X_{\ell,j}$ are defined such that:
\begin{equation}
    \mathcal{A_\ell}\left(P_{\ell,j}\begin{pmatrix} \boldsymbol{u}_\ell \\ p_\ell \end{pmatrix},\begin{pmatrix} \boldsymbol{v}_{\ell,j} \\ q_{\ell,j} \end{pmatrix}\right) = \mathcal{A_\ell}\left(\begin{pmatrix} \boldsymbol{u}_\ell \\ p_\ell \end{pmatrix},\begin{pmatrix} \boldsymbol{v}_{\ell,j} \\ q_{\ell,j} \end{pmatrix}\right), \qquad \forall \begin{pmatrix} \boldsymbol{v}_{\ell,j} \\ q_{\ell,j} \end{pmatrix} \in X_{\ell, j}.
\end{equation}

Subspace decompositions can be applied additively or multiplicatively. While the additive version is trivially parallelizable, it requires a damping factor less than 1/8 in three dimensions, which results in suboptimal convergence properties. Hence, we focus on a parallel implementation the multiplicative Schwarz smoother, defined as follows:
\begin{gather*}
\begin{aligned}
    \mathcal{P}_\ell := & I - (I - P_{\ell,J})\dotsb(I-P_{\ell,2})(I-P_{\ell,1}) \\
        & I - (I - R^T_{\ell,J}A^{-1}_{\ell,J}R_{\ell,J}\mathcal{A}_\ell)\dotsb(I-R^T_{\ell,2}A^{-1}_{\ell,2}R_{\ell,2}\mathcal{A}_\ell)(I - R^T_{\ell,1}A^{-1}_{\ell,1}R_{\ell,1}\mathcal{A}_\ell),
\end{aligned}
\end{gather*}
where $R_{\ell,j}: X_\ell \rightarrow X_{\ell,j}$ is the restriction operator and its transpose is the embedding.
\begin{figure}[tp]
\centering
    \includegraphics[width=.65\textwidth]{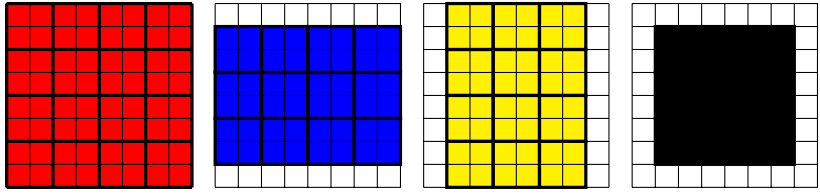}
    \caption{Non-overlapping coloring for vertex patches on regular meshes. In each color we strive to obtain a parqueting of the whole domain. Subsequent colors are obtained by shifting the first patch by one cell in each coordinate direction.}
    \label{fig:patch_coloring}
\end{figure}
To facilitate parallel execution of the multiplicative smoother, a vertex patch coloring method is applied that satisfies $\mathcal{A}_\ell$-orthogonality conditions, ensuring that $R_{\ell,i}\mathcal{A}_\ell R^T_{\ell,j}=0$. The operational stencil of $\mathcal{A}_\ell R^T_{\ell,j}$
is illustrated in Figure~\ref{fig:smoother}, where the blue areas indicate the local degrees of freedom. This allows the patches to be subdivided into non-overlapping groups as depicted in Figure~\ref{fig:patch_coloring}, which can then be processed in parallel, avoiding race conditions. Following standard literature conventions, we refer to this process of subdivision as ``colorization''. Summarizing this section, we present the colorized vertex-patch smoother in Algorithm~\ref{alg:Multiplicative_Schwarz}.
\begin{algorithm}[tp]
\caption{Colorized Vertex-Patch Smoother $S_\ell(x_\ell,b_\ell)$.}
\label{alg:Multiplicative_Schwarz}
\begin{algorithmic}[1]
\For{$c=1,\dots,n_{c}$} \Comment{Sequential loop over colors}
\For{\textbf{each } $j \in\mathcal{J}_c$} \Comment{parallel loop inside color}
\State $r_\ell \gets b_\ell - \mathcal{A}_\ell x_\ell$ \Comment{local residual}
\State $x_\ell \gets x_\ell + R^T_{\ell,j} A_{\ell,j}^{-1} R_{\ell,j} r_\ell$ \Comment{local solver}
\EndFor
\EndFor
\end{algorithmic}
\end{algorithm}

\subsection{Local solver and fast diagonalization}\label{sec:local_solver}
In computational practice, we encounter the recurring challenge of efficiently inverting local discretization matrices. Directly solving these local problems via Singular Value Decomposition (SVD) is prohibitively costly, especially for high-order finite elements. Consider the local matrix $\mathcal{A}_j$ for the $H^{\text{div}}$ interior penalty method on the subdomain $\Omega_j$ in 2D
\begin{equation}\label{eq:block_A}
    \mathcal{A}_j = \begin{pmatrix}
        A_{j} & B^T_j \\ B_j & 0
    \end{pmatrix}.
\end{equation}
With the tensor product structures of Raviart-Thomas elements, the matrix $A_{j}$ can be expressed as:
\begin{equation}
    A_{j} = \begin{pmatrix}
        L_1 \otimes M_2 + M_1 \otimes L_2 & 0 \\ 0 & L_1 \otimes M_2 + M_1 \otimes L_2
    \end{pmatrix}
\end{equation}
and the matrix $B_j$ as
\begin{equation}\label{eq:b}
   B_j = \begin{pmatrix}
        D_1 \otimes M'_2 & M'_1 \otimes D_2.
    \end{pmatrix}
\end{equation}
Here, $D_i, i=1,2$ relates to the $i$’th partial derivative arising from the divergence bilinear form $(p_h, \nabla\cdot\boldsymbol{v}_h)$, restricted to the $i$’th interval of the Cartesian product that defines the vertex patch $\Omega_j$. Similarly, $M'_i$ represent the univariate mass matrix with velocity ansatz and pressure test functions. And $L_i$, arising from $a_h(\boldsymbol{u}_h,\boldsymbol{v}_h)$ is the univariate SIPG discretization of the Laplacian, applied to the $i$’th vector component of the local velocity field, with $M_i$ being the univariate mass matrix~\cite{witte2021fast,Cui2024}. With the anisotropic tensor structure of Raviart-Thomas elements, $L_1$ and $L_2$ have different dimensions, $N\times N$ and $(N+1)\times (N+1)$, respectively. The same applies to $M_i$. Additionally, $M'$ and $D_i$ have dimensions $N \times (N+1)$ and are no longer square matrices.

Since each diagonal block of $A_{j}$ admits the separable Kronecker representation, its inverse can be computed efficiently using fast diagonalization methods as proposed by ~\cite{lynch1964direct}
\begin{equation}
\label{eq:inverse2d}
    A_{j}^{-1} = S_2 \otimes S_1 (\Lambda_2 \otimes I + I \otimes \Lambda_1)^{-1} S_2^T \otimes S_1^T
\end{equation}
where $S_i, i=1,2$ is the matrix of eigenvectors to the generalized eigenvalue problem in the given tensor direction $i$:
\begin{equation}
\label{eq:fast_inverse}
    L_iS_i=M_iS_i\Lambda_i , \quad i=1, 2,
\end{equation}
and $\Lambda_i$ is the diagonal matrix representing the generalized eigenvalues $\lambda_i$. 

While local operators $\mathcal{A}_j$ do not lend themselves to direct fast inversions, we employ the Schur complement for the solution:
\begin{gather*}
    \begin{aligned}
        B_j^TA_{j}^{-1}B_j P &= B^TA_{j}^{-1}F - G, \\
        A_{j}U &= F - BP,
    \end{aligned}
\end{gather*}
where $(F, G)^T$ represents the right-hand side vectors. 
For the pressure component $P$, we solve the first equation $SP =B^TA_{j}^{-1}F - G $ using the Conjugate Gradient method, where $S$ is the Schur complement with the formula $S = B_j^TA_{j}^{-1}B_j$. During the iterations, there is no need to explicitly compute $S$. Instead, we compute $Sv = B_j^T\left(A_{j}^{-1}\left(B_jv\right)\right)$, where $B$ has a low-rank representation as shown in \eqref{eq:b}. Meanwhile, the operation of $A_j^{-1}$ can be performed efficiently using the fast diagonalization method~\cite{lynch1964direct} as in~\cite{Cui2024}. Finally, we apply the fast diagonalization method again to solve the second equation for the velocity component $U$.

\subsection{Efficiency and convergence}\label{sec:tests}

The motivation for higher-order methods lies in their computational intensity, which allow for achieving higher accuracy given a fixed budget of degrees of freedom, or more commonly, attaining a certain level of accuracy with fewer degrees of freedom. Our ultimate goal is to minimize computation time while achieving the prescribed accuracy.
Our numerical evidence indicates that the proposed smoother produces a highly efficient multigrid method in terms of iteration count. Summarizing our findings, we observe that the iteration count is independent of the refinement level $L$. Furthermore, our results demonstrate that the convergence steps decrease progressively with the increase in polynomial degree.
 
We solve the Stokes problem~\eqref{eq:model} with the GMRES method preconditioned by a multigrid V-cycle varying finest level $L$ to a relative accuracy of $10^{-8}$. The right-hand side $\boldsymbol{f}$ in the Stokes equations ~\eqref{eq:model} is manufactured from reference solutions $\boldsymbol{u}$ and $p$,
\begin{gather}
    \begin{aligned}
        \boldsymbol{u}(x,y) =& \begin{pmatrix} \partial_y\psi \\ -\partial_x \psi \end{pmatrix}, \\
        p(x,y) =& \cos(2\pi x) \cos(2\pi y)
    \end{aligned}
\end{gather}
in two dimensions and 
\begin{gather}
    \begin{aligned}
        \boldsymbol{u}(x,y,z) =& \begin{pmatrix}
        \partial_y\psi + \partial_z\psi \\ -\partial_x\psi -\partial_z\psi \\ -\partial_x\psi + \partial_y \psi
    \end{pmatrix}, \\
        p(x,y,z) =& \cos(2\pi x) \cos(2\pi y) \cos(2\pi z)
    \end{aligned}
\end{gather}
in three dimensions, with 
\begin{gather}
    \begin{aligned}
        &\psi(x,y) = \varphi(x) \varphi(y) &\quad\text{in 2D}, \\
        &\psi(x,y,z) = \varphi(x) \varphi(y) \varphi(z) &\quad\text{in 3D},
    \end{aligned}
\end{gather}
where
\begin{equation*}
    \varphi(x) = \frac{x^2(x-1)^2}{\sqrt{2\pi\sigma^2}}\exp{\left(-\frac{(x-\mu)^2}{\sigma^2}\right)}.
\end{equation*}
We measure efficiency of the preconditioners by reduction of the Euclidean norm of the residual $\|r_n\|$ after $n$ steps compared to the initial residual norm. Since the (integer) number of iteration steps depends strongly on the chosen stopping criterion,
we define the \emph{fractional iteration count} $\nu$ by
\begin{equation*}
    \nu = -8 \log_{10} \bar{r}, \quad \bar{r} = \left(\frac{\|r_n\|}{\|r_0\|} \right)^{1/n}.
\end{equation*}

Since the pressure component is solved iteratively with on average 15 CG iterations, the local solver $A_j$ within the vertex-patch smoother with Schur complement approach operates as an approximate solver.
Given the inherent approximations in this approach, the local solver with the pseudo-inverse of $A_j$ (omitting the smallest singular value, which is zero) serves as a benchmark (referred to as \emph{Direct solver}). This direct solver configuration fixes the mean value of the pressure implicitly, providing a reliable reference for assessing the numerical efficiency of the Schur complement method.

\begin{table}[tp]
\centering
\caption{Fractional iteration count of GMRES iterations for Stokes equation in two dimensions, GMRES solver with relative accuracy $10^{-8}$ preconditioned by multigrid.}
\begin{tabular}{c|cccccc|cccccc}
\hline
& \multicolumn{6}{c|}{Direct solver} & \multicolumn{6}{c}{Schur complement} \\ 
\hline
$L$ & $\mathbb{Q}_{3}$ & $\mathbb{Q}_{4}$ & $\mathbb{Q}_{5}$ & $\mathbb{Q}_{6}$ & $\mathbb{Q}_{7}$ & $\mathbb{Q}_{8}$ & 
    $\mathbb{Q}_{3}$ & $\mathbb{Q}_{4}$ & $\mathbb{Q}_{5}$ & $\mathbb{Q}_{6}$ & $\mathbb{Q}_{7}$ & $\mathbb{Q}_{8}$\\
\hline
    3 & 2.8 & 2.4 & 2.1 & 2.2 & 1.9 & 2.1 & 2.6 & 2.3 & 2.0 & 2.2 & 2.0 & 2.0\\
    4 & 3.2 & 2.4 & 2.3 & 2.2 & 2.1 & 2.0 & 3.2 & 2.5 & 2.4 & 2.3 & 2.1 & 2.2\\
    5 & 3.2 & 2.3 & 2.3 & 2.0 & 2.0 & 1.9 & 3.3 & 2.5 & 2.5 & 2.4 & 2.0 & 2.2\\
    6 & 3.0 & 2.2 & 2.2 & 1.9 & 1.9 & 1.8 & 3.3 & 2.5 & 2.5 & 2.3 & 1.9 & 2.0\\
    7 & 2.9 & 1.9 & 1.9 & 1.8 & 1.8 & 1.7 & 3.3 & 2.5 & 2.4 & 2.3 & 1.9 & 2.0\\
    8 & 2.8 & 1.9 & 1.8 & 1.8 & 1.8 & 1.8 & 3.0 & 2.5 & 2.4 & 2.3 & 1.9 & 1.9 \\
    9 & 2.7 & 1.9 & 1.9 & 1.8 & 1.8 & 1.8 & 2.9 & 2.4 & 2.4 & 2.4 & 1.9 & 1.9 \\
\hline
\end{tabular}
\label{tab:exact-2d}
\end{table}

\begin{figure}[tp]
\centering
\includegraphics[width=.45\textwidth]{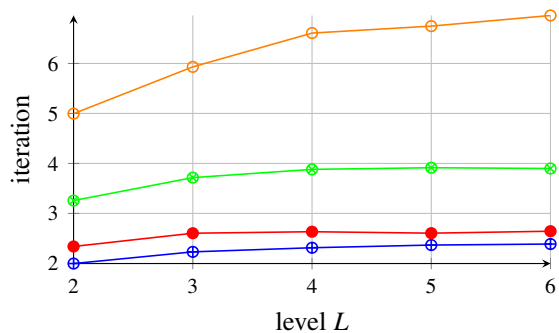}
\caption{Fractional iteration count of GMRES iterations for Stokes equation in three dimensions, GMRES solver with relative accuracy $10^{-8}$ preconditioned by multigrid using Schur complement local solver, polynomial degree $k=1$(\ref{pgfplots:p1}),$2$(\ref{pgfplots:p2}),$3$(\ref{pgfplots:p3}),$4$(\ref{pgfplots:p4})}
\label{fig:frac_its}
\end{figure}

Convergence results in two dimensions are presented in Table~\ref{tab:exact-2d}. The data shows that the vertex-patch smoother facilitates uniform convergence in both cases. Notably, the number of iterations required for convergence does not depend on the polynomial degree $k$ or the mesh refinement level $L$, and it even decreases with an increase in polynomial degree. This behavior underscores the solver's robustness and efficiency, classifying it as highly effective for numerical simulations.
Further validation is conducted in three dimensions in Figure~\ref{fig:frac_its}, where we use the same solver settings and stopping criteria. Consistent with the 2D outcomes, the number of iterations remains independent of the mesh refinement level. From polynomial degree 3 onward, the solver requires no more than three steps to converge, approaching the efficiency of a direct solver.
\begin{figure}[tp]
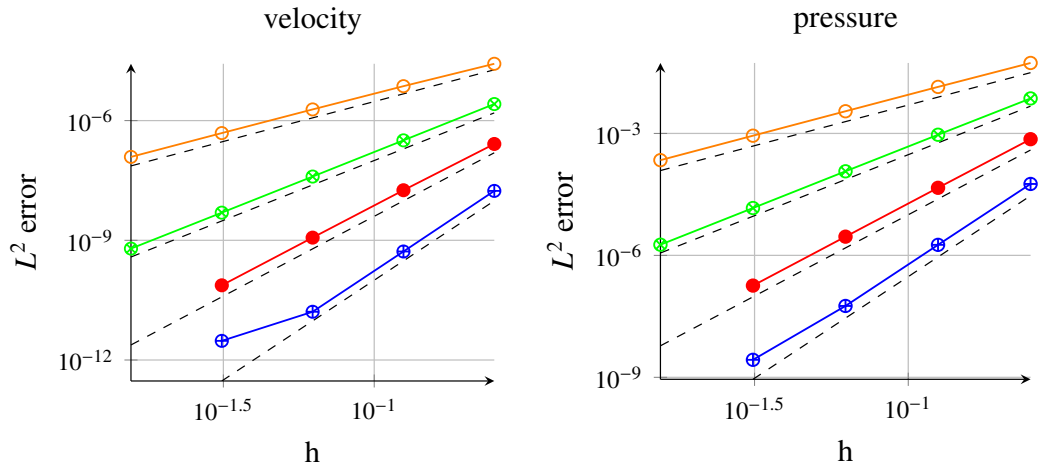

\centering
\includegraphics[width=.4\textwidth]{Figures/convergence3d.tikz}
\hspace{0.3cm}
\includegraphics[width=.4\textwidth]{Figures/convergence3dp.tikz}
\caption{Convergence of 3D Stokes problem, $\text{RT}_k$ element with polynomial degree $k=1$(\ref{pgfplots:p1}),$2$(\ref{pgfplots:p2}),$3$(\ref{pgfplots:p3}),$4$(\ref{pgfplots:p4}).}
\label{fig:convergence3d}
\end{figure}
\begin{figure}[tp]
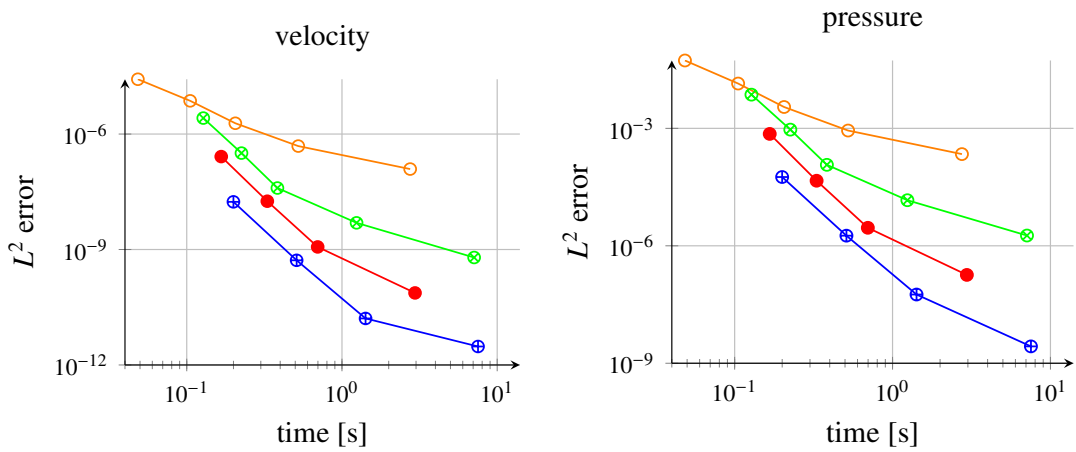

\centering
\includegraphics[width=.42\textwidth]{Figures/convergence_time.tikz}
\hspace{0.2cm}
\includegraphics[width=.42\textwidth]{Figures/convergence_timep.tikz}
\caption{$L_2$ error vs. total run time of 3D Stokes problem, $\text{RT}_k$ element with polynomial degree $k=1$(\ref{pgfplots:p1}),$2$(\ref{pgfplots:p2}),$3$(\ref{pgfplots:p3}),$4$(\ref{pgfplots:p4}).}
\label{fig:convergence3d_time}
\end{figure}

Figure~\ref{fig:convergence3d} illustrates the optimal convergence orders for both velocity and pressure components for $\text{RT}_k, k = 1, 2, 3$ and 4. Given our objective to minimize computation time, we further examine the comparison between error and solution time as shown in Figure~\ref{fig:convergence3d_time}. The results indicate a significant advantage of higher-order methods in terms of runtime. For a given target accuracy, higher-order methods can substantially reduce the computation time.

\section{GPU Implementation}\label{sec:gpu}

Building on our previous work~\cite{Cui2023,Cui2024}, which demonstrated the efficacy of our algorithms in addressing the Poisson problem, we continue to evolve these design principles to enhance our current study problem. We use implementations based on the \texttt{deal.II}~\cite{dealII95} finite element library and in particular its \texttt{Multigrid} framework. All numerical experiments in this work are performed on a single NVIDIA Ampere A100 SXM4 GPU with 80GB of high-speed HBM2e memory for VRAM which provides 2 TB/s peak memory bandwidth and 8.7 TFLOPS/s peak double precision performance at 1.27 GHz, hosted on a system with two AMD EPYC 7282 16-Core processors.

First, a crucial aspect of optimizing performance is minimizing memory transfers between the host and device. To this end, the entire solver has been implemented directly on the GPU. During the setup phase, all necessary vectors and data structures are allocated and initialized on the GPU. This setup ensures that there is no need for data communication between the CPU and GPU during the solution of the linear system, enhancing performance and reducing latency.
Then, achieving high performance on the GPU also depends on effective utilization of on-chip memory, commonly referred to as \emph{Shared Memory}. Our approach adheres to a typical CUDA programming pattern: 1) Data is loaded from the Video Random Access Memory (VRAM) to shared memory, 2) Computations are performed on the data within shared memory, and 3) Results are stored back to VRAM from shared memory. Proper synchronization ensures that data integrity is maintained throughout these steps. 

\begin{figure}[tp]
\centering
\includegraphics[width=.3\textwidth]{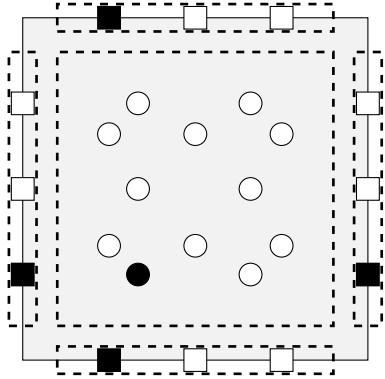}
\caption{Degrees of freedom for $RT_2$ in 2D visualized as support points with stored indices in black.}
\label{fig:index_storage}
\end{figure}

Given the limited cache and memory capacities of modern GPUs compared to CPUs, it's important to use an efficient index storage scheme. This helps maximize the use of shared memory and manage larger problem sizes. Our approach stores only the first index of each entity within a cell, using it to directly load vector entries into shared memory. Figure~\ref{fig:index_storage} illustrates this index storage method for a 2D Raviart-Thomas element of order 2.

In multigrid methods, operators are typically implemented by iterating over all cells (or patches in vertex-patch smoothers). The operations performed on each cell are independent of those on other cells, allowing for parallel execution. Within each cell, all degrees of freedom are evaluated simultaneously, eliminating the need for inter-cell communication. This computational pattern is well-suited to the CUDA programming model, where each CUDA thread corresponds to a degree of freedom, and each CUDA thread block corresponds to an cell. This setup facilitates two levels of parallelism, thereby enhancing computational efficiency.

Our focus is on the kernel performance of the most computationally demanding tasks, namely the matrix-free evaluation of the Stokes operator and the local solvers within the smoothing operator. In the subsequent sections, we detail the continuous optimizations undertaken to achieve better performance. We compare the TFLOPS/s achieved by each kernel version with the empirically determined roofline performance model. Additionally, to provide a clearer understanding of the algorithm's performance, we evaluate the throughput of operations in terms of degrees of freedom per second (DoF/s).

\subsection{Stokes operator kernels}\label{sec:stokes_kernels}

\begin{figure}[tp]
\centering
\includegraphics[width=.45\textwidth]{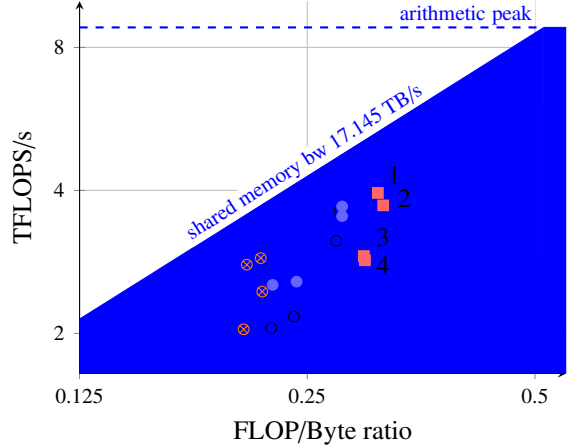}
\caption{Shared memory roofline performance model for Stokes operator with various GPU implementations in 3D. Comparison of different kernels: Padding (\ref{pgfplots:padding_marks}), Basic (\ref{pgfplots:basic_marks}), Shared memory (\ref{pgfplots:cache_marks}) and Conflict-free (\ref{pgfplots:cf_marks}) for $\text{RT}_k$ element with polynomial degree $k=1,2,3,4$ on NVIDIA A100 GPU at 1.27 GHz.}
\label{fig:roofline}
\end{figure}

The implementation of the Stokes operator is based on a cell-wise approach. Since the face terms require information from neighboring cells, the data is loaded in a patch-wise manner. In this way, we compute the contribution for each cell as well as part of the face integrals, ensuring that each face in the mesh is computed once and only once. Input vectors are loaded into shared memory to ensure that local operations avoid the random long latency global memory access. The results shown in Figure~\ref{fig:roofline} help illustrate the impact of each optimization on the overall performance of the Stokes operator. Below, we provide details for each kernel, along with any optimizations applied within that kernel. Due to the anisotropic nature of Raviart–Thomas elements, our matrices are no longer square, necessitating different approaches for implementing local operations:

\textbf{\textit{Padding Kernel}}. In this kernel, the constant data, such as 1D objects and mapping data, are stored in constant memory. We pad the matrix to restore its square structure for computation. Specifically, we only need zero-padding in the orthogonal direction. This method is straightforward and easy to implement. However, it has its drawbacks. Although padding is only applied in the orthogonal directions, see Figure~\ref{fig:RT2}, it still has a significant impact in three-dimensional cases, exacerbating the already substantial pressure on shared memory. Additionally, the extra operations introduced by the larger-dimensional matrix computations also affect the overall performance. Consequently, this kernel achieves only 2.88 TFLOPS/s.

\textbf{\textit{Basic Kernel}}. In this kernel, we match the number of threads to the matrix sizes, with one thread handling one ``column'' degree of freedom in three dimensions. This requires careful attention to the corresponding matrix operations for different directions. We can still use template parameters on loop bounds, allowing the compiler to fully unroll loops and rearrange operations to improve instruction flow. This kernel achieves up to 3.62 TFLOPS, a 25\% improvement over \emph{Padding} kernel.

\textbf{\textit{Shared Memory Kernel}}. Unlike the previous kernels, this kernel loads all data, including constants, into shared memory. While accessing constant memory is faster than shared memory, constant memory access has stricter limitations: optimal performance is only achieved when all threads in a warp access the same address. Therefore, shared memory is a more suitable choice, and the performance gain stems from this approach.

\textbf{\textit{Conflict-Free Kernel}}: The final challenge in shared memory is bank conflicts. Building on our previous work~\cite{Cui2023}, by rearranging the order of tensor contraction operations and data layout in shared memory, our implementation ensures an optimal access pattern. This leads to significant performance improvement, reaching 3.95 TFLOPS/s, which is 45\% of the double-precision peak performance.

\textbf{Roofline analysis.} The traditional roofline model assumes VRAM bandwidth as the limiting factor. However, due to the extensive use of shared memory, we assume that shared memory bandwidth is the true limiting factor. Inspired by~\cite{swirydowicz2019acceleration,Cui2023,cui2024acceleration}, we consider a roofline model that specifically focuses on shared memory bandwidth to more accurately understand the factors affecting performance. The bandwidth of shared memory can be estimated using the formula:
\begin{equation*}
    B = \# \text{SMs} \times \# \text{banks} \times \text{word length} \times \text{clock speed}.
\end{equation*}
For NVIDIA A100 GPU, the corresponding bandwidth is $108\cdot 32 \cdot 4 \cdot 1.27=17.145 \text{TB/s}$, which aligns with the measured shared bandwidth of 127.7 bytes/clk/SM as reported in~\cite{sun2022dissecting}. 

First, from Figure~\ref{fig:roofline}, we can see that the \emph{Padding} kernel appears on the most left, exhibiting the lowest FLOP/Byte ratio due to the additional zero-padding introduced. The \emph{Basic} and \emph{Shared Memory} kernels have almost identical arithmetic intensity; however, their performance is reduced because of the presence of bank conflicts. Each bank conflict forces an additional new memory transaction. The more transactions that are required, the more unused words are transferred alongside the words accessed by the threads, also reducing the instruction throughput.
Finally, the \emph{Conflict-free} kernel significantly enhances performance, achieving 40\% of peak performance for degree $k = 1, 2$. This is lower than the typical 90\% achieved by large dense matrix-matrix multiplication \footnote{By using the vendor library cuBlas~\cite{cublas}, we were able to achieve over 95\% of the theoretical peak performance on a single NVIDIA A100 GPU.}. The primary reasons for this lower performance are the unstructured read/write access to the global vector and the operations at quadrature points, which introduce non-floating point instructions and affect the instruction pipeline. For higher-order elements, performance noticeably declines due to excessive use of shared memory, impacting achieved occupancy. Only one thread block can be assigned to a streaming multiprocessor, leading to a drop in performance. This is verified by NVIDIA's profiling tool \texttt{Nsight Compute}~\cite{cudaNcu}.

\begin{figure}[tp]
\centering
\includegraphics[width=.45\textwidth]{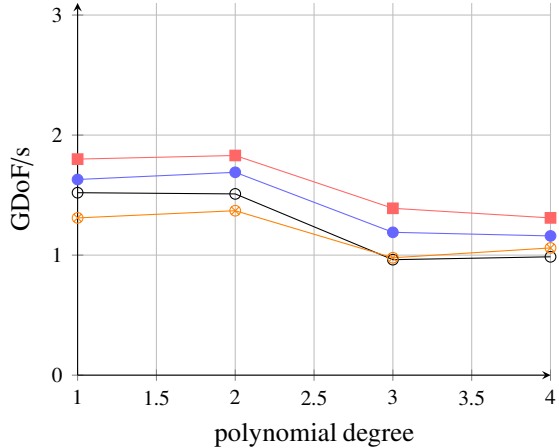}
\caption{Throughput measured in GDoF/s for Stokes operator with various implementations in 3D: Padding (\ref{pgfplots:padding}), Basic (\ref{pgfplots:basic}), Shared memory (\ref{pgfplots:cache}) and Conflict-free (\ref{pgfplots:cf}).}
\label{fig:perf_Ax}
\end{figure}

TFLOPS/s does not fully convey the actual speed of operator evaluation. To better understand algorithm performance, we also evaluate the throughput in terms of degrees of freedom per second (DoF/s). As shown in Figure~\ref{fig:perf_Ax}, the improvements in TFLOPS/s are also reflected in DoF/s, exhibiting similar performance patterns, with a peak throughput approaching 2 GDoF/s in the best case. We also observe a performance decline for $k = 3, 4$.
Despite the challenge with shared memory usage, the performance improvements in TFLOPS/s translate directly into increased throughput in DoF/s, demonstrating the effectiveness of the optimizations across different metrics.

\subsection{Local solver kernels}
\begin{figure}[tp]
\centering
\includegraphics[width=.45\textwidth]{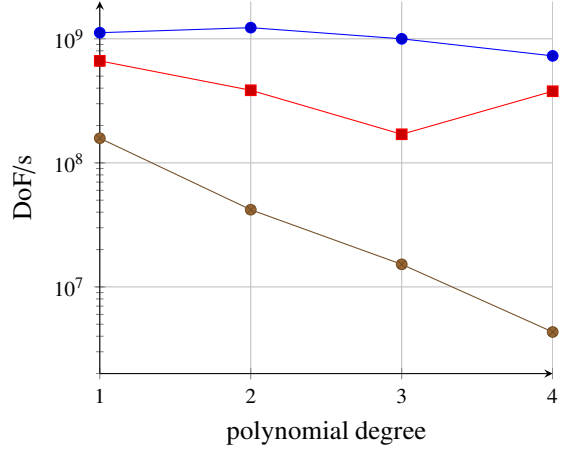}
\caption{Throughput measured in DoF/s for local solver per color with various implementations in 3D: Conflict-free (\ref{pgfplots:tp_cf}), Schur complement (\ref{pgfplots:tp}) and Direct (\ref{pgfplots:direct}). The problem sizes are
between 28 and 131 million DoFs.}
\label{fig:perf_smoother}
\end{figure}

Another crucial operation is the local solver in the smoothing operator. To avoid explicit inversion and multiplication by the inverse, which involves $\mathcal{O}(k^{2d})$ operations, we use the Schur complement as described in Section~\ref{sec:local_solver}. Additionally, in three-dimensional cases, the memory required for the local inverse matrix easily exceeds the 164KB capacity of shared memory. As shown in Figure~\ref{fig:perf_smoother}, accessing global memory significantly impacts performance.

By employing the Schur complement in combination with fast diagonalization, all operations can be completed within shared memory, resulting in relatively stable performance. Similarly, we can apply the conflict-free kernel techniques discussed in Section~\ref{sec:stokes_kernels} to eliminate bank conflicts, as the fundamental operations are tensor products.
Comparing the two Schur Complement kernels (Conflict-free and Schur Complement) in Figure~\ref{fig:perf_smoother}, eliminating bank conflicts results in at least a 1.5-fold performance improvement, with nearly a sixfold increase at $\text{RT}_3$. Unlike the Stokes operator, the conflict-free access pattern is particularly crucial for the Schur Complement kernel due to the extensive tensor product operations involved in solving the pressure component using iterative methods. Below, we analyze the causes of bank conflicts and why they are most pronounced at $\text{RT}_3$, using polynomial degree $k=1,3$ as examples.

\begin{figure}
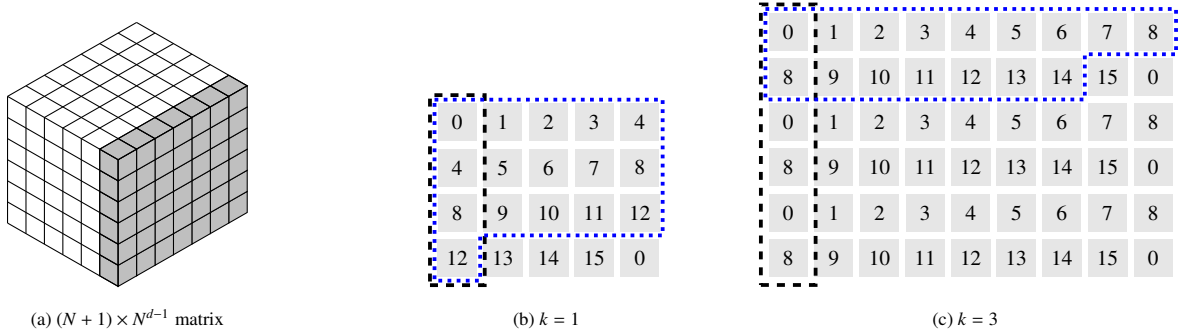

  \begin{subfigure}[b]{0.33\textwidth}
\centering
\includegraphics{Figures/mat_d2.tikz}
\caption{$(N+1)\times N^{d-1}$ matrix}
  \end{subfigure}
  \begin{subfigure}[b]{0.33\textwidth}
\centering
\includegraphics{Figures/bank_d2_k1.tikz}
\caption{$k=1$}
  \end{subfigure}
  \begin{subfigure}[b]{0.33\textwidth}
\centering
\includegraphics{Figures/bank_d2.tikz}
\caption{$k=3$}
  \end{subfigure}
\caption{Data layout and shared memory access pattern for tensor contraction with direction 2.}
\label{fig:bank_conflict}
\end{figure}

\begin{figure}[tp]
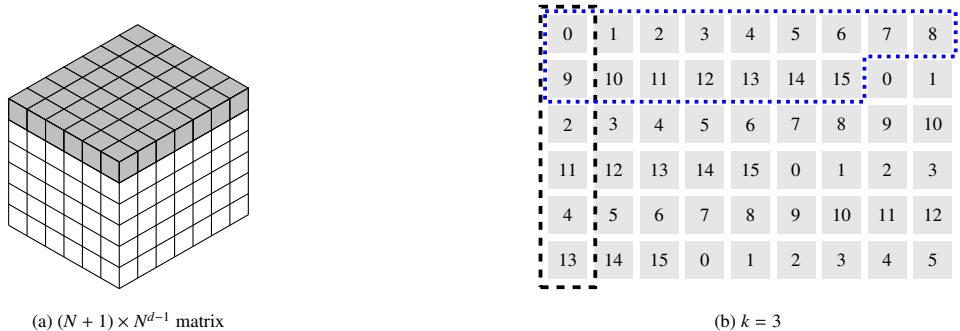

  \begin{subfigure}[b]{0.49\textwidth}
\centering
\includegraphics{Figures/mat_d1.tikz}
\caption{$(N+1)\times N^{d-1}$ matrix}
\end{subfigure}
  \begin{subfigure}[b]{0.49\textwidth}
\centering
\includegraphics{Figures/bank_d1.tikz}
\caption{$k=3$}
\end{subfigure}
\caption{Conflict-free shared memory access pattern for tensor contraction with direction 1.}
\label{fig:tensor_mul}
\end{figure}

In Figure~\ref{fig:bank_conflict}(a), we depict the reshaped input vector corresponding to an $(N+1)\times N^{d-1}$ matrix, considering a row-major storage format, where elements within each row are stored contiguously in memory. The operator evaluation involves different directions, implying matrix multiplications in various forms. The figure shows a matrix-matrix multiplication when operating in direction 2.
Figures~\ref{fig:bank_conflict}(b) and (c) illustrate the bank numbers corresponding to the matrix elements in shared memory for $k=1$ and $k=3$ respectively, considering double precision, which results in 16 banks (32 banks for single precision).

For $\text{RT}_3$, There are no bank conflicts during local computation, as shown by the access pattern within the black dashed box. Bank conflicts only occur when storing the computed results back to shared memory (blue dashed box).
For $\text{RT}_3$, When performing matrix multiplication, all threads accessing the first column result in a four-way bank conflict, meaning these accesses are serialized, reducing efficiency to 25\% of the conflict-free access rate. Additionally, storing the computed results back to shared memory also involves bank conflicts. This explains why $\text{RT}_3$ is most severely affected and sees the most significant performance improvement when bank conflicts are eliminated.

\subsection{Mixed Precision}

When solving linear systems of equations, whether by iterative or direct methods, maintaining high floating-point precision is crucial for obtaining accurate results, as rounding errors can accumulate during computation. Previous study~\cite{goddeke2007performance,Cui2023} compared double-precision algorithms and mixed-precision algorithms for solving linear systems using Conjugate Gradient (CG) and multigrid methods with a Jacobi smoother. The proposed mixed-precision method performed exceptionally well, maintaining the same accuracy as the reference solver running in double precision.

\begin{figure}[tp]
\centering
\includegraphics[width=.75\textwidth]{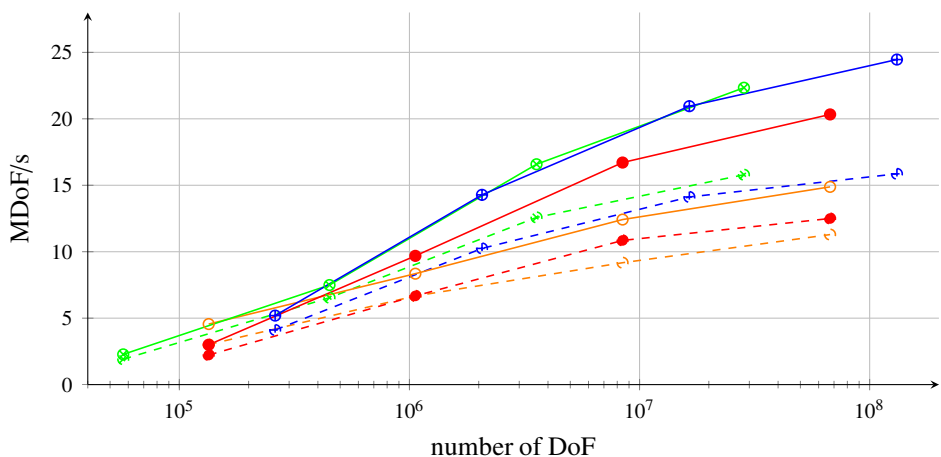}
\caption{Comparison of Double (dashed) and Mixed (solid) Precision performance for solving the Stokes equation with preconditioned GMRES method in three dimensions, $\text{RT}_k$ element with polynomial degree $k=1$(\ref{pgfplots:p1}),$2$(\ref{pgfplots:p2}),$3$(\ref{pgfplots:p3}),$4$(\ref{pgfplots:p4}).}
\label{fig:throughput}
\end{figure}
In our mixed precision approach, the multigrid V-cycle is fully done in single precision and the outer GMRES iteration is done in double precision. The format is converted when entering and exiting the V-cycle. Figure~\ref{fig:throughput} shows the performance results for solving the Stokes equations, measured in millions of degrees of freedom per second. For small-scale problems, both double-precision and mixed-precision solutions exhibit relatively low performance, as the problem's parallelism is insufficient to fully saturate the hardware. For higher-order elements, the mixed-precision approach improves solution speed by a factor of 1.6. The performance curves rise smoothly, thanks to the optimized solver and the flexibility of CUDA, which ensures that performance is not limited by CPU cache sizes. Additionally, despite the performance decrease of the Stokes operator with higher-order elements, as shown in Figure~\ref{fig:perf_Ax}, the reduction in the number of iterations (in Figure~\ref{fig:frac_its}) significantly shortens the solution time, achieving performance exceeding $2 \cdot 10^7$ degrees of freedom per second. 

\subsection{Comparison against Poisson Problem}


In our final evaluation of the efficiency of the GPU implementation in this work, we compare the performance of solving the Stokes problem with that of the 3D Poisson problem. There, we have proposed a highly optimized implementation in~\cite{Cui2024}.
The Poisson problem is discretized using the symmetric interior penalty method, with a GMRES solver preconditioned by a multigrid method incorporating a vertex-patch smoother. The operator application uses sum factorization and the local solvers are optimized by using fast diagonalization.

We actually compare to the results of two different versions of the patch smoother for the Poisson problem: the first one is organized like the smoother in this work, alternating between a global residual computation and application of the local solvers in a whole set of nonoverlapping patches. The \emph{fused} kernel computes the residual on the fly before each local solve, which reduces operations and global data transfers. However, in order to do so, the local solver uses a subspace with zero boundary conditions and ignores parts of the residual. Compared to the original version, the resulting inexact local residuals lead to a higher number of iterations for convergence, as illustrated in Table~\ref{tab:vs_poisson}. Nevertheless, the \emph{fused} smoother improves the overall performance of the solver.
We cannot apply this optimization to the Stokes problem, as experiments show it results in extremely slow or even stalled convergence.

The stopping criterion for both problems is a relative residual reduction of $10^{-6}$. For the Poisson problem, $\mathbb{Q}_k$ element with polynomial degrees $k=4,7$ are chosen, while for the Stokes problem, $\text{RT}_2$ and $\text{RT}_4$ Raviart-Thomas elements are selected, ensuring a similar workload per patch.

\begin{table}[tp]
\caption{Comparison against solving Poisson problem in 3D with same level of refinement, i.e. on $32^3$ mesh.}
\centering
\begin{tabular}{lcccccc}
\toprule
& \multicolumn{2}{c}{Poisson} & \multicolumn{2}{c}{Poisson (Fused)} & \multicolumn{2}{c}{Stokes} \\
  \cmidrule(lr){2-3} \cmidrule(lr){4-5} \cmidrule(lr){6-7} 
  & $\mathbb{Q}_4$ & $\mathbb{Q}_7$ & $\mathbb{Q}_4$ & $\mathbb{Q}_7$ & $\text{RT}_2$ & $\text{RT}_4$ \\
\midrule
\text{DoF [per patch]} & 1000 & 4096 & 1000 & 4096 & 972 & 4300 \\
\text{Mat-vec [ns/DoF]} & 0.374 & 0.676 & 0.374 & 0.676 & 0.535 & 0.813 \\
\text{Smoothing [ns/DoF]} & 3.968 & 6.944 & 1.068 & 1.866 & 7.576 & 10.582 \\
\text{Iterations} & 2 & 2 & 4 & 5 & 3 & 3 \\
\text{Time to solution [ns/DoF]} & 28.875 & 49.751 & 15.060 & 31.140 & 63.694 & 79.365 \\
\bottomrule
\end{tabular}
 \label{tab:vs_poisson}
\end{table}
In Table~\ref{tab:vs_poisson}, we make two kinds of comparison. First, we compare the Poisson and Stokes solvers in a setting where the number of degrees of freedom on each patch and hence the global number of degrees of freedom are roughly equal. Hence, we can assess the drop of performance between the separable structure of the Laplacian and the much more complex Stokes operator. Comparing $\mathbb{Q}_4$ with $\text{RT}_2$ and $\mathbb{Q}_7$ with $\text{RT}_4$, respectively, we see that matrix-free operator evaluations, take 1.2 to 1.4 times longer for the Stokes problem. This performance drop can be attributed to the anisotropic tensor structure of the Raviart-Thomas elements. Additionally, the operator $\mathcal{A}$ in~\eqref{eq:block_A} includes several sub-blocks, leading to reduced branch efficiency. Nevertheless, this result indicates almost the same efficiency.

Comparing the performance for a single smoothing step with the more similar {nonfused} Poisson kernel, the Stokes smoother performs surprisingly well with timing increased by less than a factor two in spite of the iterative local solver. Comparing to the fused kernel, the performance drop in each step is significantly higher. Nevertheless, in both cases this translates to a total time to solution which is also roughly twice that of the Poisson problem with similar number of degrees of freedom. As pointed out above, we currently do not have a fused kernel available for the Stokes problem.

The other comparison is between a Stokes and a Poisson solver where the discretization of the elliptic operator is of the same size, since this allows a comparison to a pressure correction or a block preconditioner scheme. Due to the anisotropy of the Raviart-Thomas element, this is not completely feasibly, but the comparison of the $\mathbb{Q}_4$ element with $\text{RT}_4$ seems fair. Here, the fused Poisson kernel manages to produce a solution in 15~nsec/DoF, while the Stokes kernel requires about 80~nsec/DoF. This is roughly a factor five. If we compare to a pressure correction scheme where three Poisson problems for the velocity and one for the pressure must be solved, we are only slightly slower, but do not produce the perturbations of the pressure correction.


\section{Conclusions}\label{sec:con}
In this paper, we have described a vertex-patch smoother for solving Stokes problems with a matrix-free implementation that efficiently utilizes the high performance and memory bandwidth of modern graphics processing units. By exploiting the tensor product structure of the Raviart–Thomas elements, the resulting finite element operators are applied using sum-factorization, and obtain throughput of over one billion degrees of freedom per second on a single NVIDIA A100 GPU. Moreover, efficient local solver through fast diagonalization minimizes computational and space complexity, resulting in a significant gain in algorithm efficiency. The robustness of preconditioners with respect to both refinement level and polynomial degree was demonstrated in both two and three dimensions, showcasing not only uniform, but also remarkably fast convergence speed, thus verifying the accuracy and efficiency of the method.




\section*{Acknowledgements}
Cu Cui was supported by the China Scholarship Council (CSC) under grant no. 202106380059.
Guido Kanschat was supported by the Deutsche Forschungsgemeinschaft (DFG, German Research Foundation) under Germany's Excellence Strategy EXC 2181/1 - 390900948 (the Heidelberg STRUCTURES Excellence Cluster)

\bibliographystyle{apalike}
\bibliography{references}

\end{document}